\documentclass[1p]{elsarticle}
\usepackage{amsfonts}
\usepackage{amsmath}
\usepackage{amssymb}
\newtheorem{theorem}{Theorem}

\newtheorem{cor}{Corollary}

\newtheorem{remark}{Remark}

\def\qed{\hfill $\Box$}

\journal{Statistics and Probability Letters}

\begin{document}

\begin{frontmatter}%

\title{Some Refinements of Large Deviation Tail Probabilities
}
\author[gy]{L\'aszl\'o Gy\"orfi}
\ead{gyorfi@cs.bme.hu}
\author[ha]{Peter Harremo\"es\corref{cor1}}
\ead{harremoes@ieee.org}
\author[tu]{G\'abor Tusn\'ady}
\ead{tusnady.gabor@renyi.mta.hu}

\cortext[cor1]{Corresponding author}

\address[gy]{
Budapest University of Technology and Economics,
Budapest, Hungary}
\address[ha]{Copenhagen Business College, Copenhagen, Denmark }
\address[tu]{R\'{e}nyi Institute of Mathematics, Budapest, Hungary }

\begin{abstract}
We study tail probabilities via some Gaussian approximations. Our
results make  refinements to large deviation theory. The proof
builds on classical results by Bahadur and Rao. Binomial
distributions and
their tail probabilities are discussed in more detail.%
\end{abstract}%


\end{frontmatter}%

\section{Introduction}

Let $X_{1},\dots,X_{n}$ be i.i.d. random variables such that the
moment generating function $\mathbf{E}\left[  \exp\left(  \beta
X_{1}\right) \right]  $ is finite in a neighborhood of the origin.
For fixed $\mu>\mathbf{E}\left[  X_1\right]$, the aim of this
paper is to approximate the tail distribution:
\[
P_{n,\mu}:=\mathbf{P}\left\{
\frac{1}{n}\sum_{i=1}^{n}X_{i}\geq\mu\right\} .
\]

If $\mu$ is close to the mean of $X_{1}$ one would usually
approximate $P_{n,\mu}$ by a tail probability of a Gaussian random
variable. If $\mu$ is far from the mean of $X_{1}$ the tail
probability can be estimated using large deviation theory.
According to the Sanov theorem the probability that the deviation
from the mean is as large as $\mu$ is of the order $\exp\left(
-nD\right)  $ where $D$ is a constant.
Bahadur and Rao \cite{Bahadur1960a} improved the estimate of this
large deviation probability, and the goal of this paper is to
extend the Gaussian tail approximations into situations where one
normally uses large deviation techniques.

Let $\phi$ and $\Phi$ be the \emph{density function and the distribution function of the standard Gaussian}, respectively.
Let $P_{0}$ denote a probability measure describing the
distribution of a random variable $X.$ Consider the 1-dimensional
\emph{exponential family}
$\left(  P_{\beta}\right)  $ based on $P_{0}$ and given by%
\[
\frac{\mbox{\textrm{d}}P_{\beta}}{\mbox{\textrm{d}}P_{0}}\left(
x\right)
=\frac{\exp\left(  \beta\cdot x\right)  }{\mathcal{Z}\left(  \beta\right)  }%
\]
where the denominator is the \emph{moment generating function}
(partition
function) given by%
\[
\mathcal{Z}\left(  \beta\right)  =\int\exp\left(  \beta\cdot
x\right) ~\mbox{\textrm{d}}P_{0}x =\mathbf{E}\left[
\mathrm{e}^{\beta  X}\right].
\]
The mean value of $P_{\beta}$ is
\begin{equation}
\label{**}
 \frac{\mathcal{Z}^{\prime}\left(  \beta\right)
}{\mathcal{Z}\left(
\beta\right)  }%
\end{equation}
and the range of this function will be denoted $M$ and will be
called \emph{the mean value range} of the exponential family.

For $\mu$ in interior of $M$ the \emph{maximum likelihood
estimate} $\hat{\beta}\left(  \mu\right) $ equals the $\beta$ such
that the mean value of $P_{\beta}$ equals $\mu$,
which in this case is the average of the i.i.d. samples. Put
$P^{\mu}=P_{\hat{\beta}\left(  \mu\right)  }.$ An equivalent
definition of $\hat{\beta}\left(  \mu\right) $ can be as the
solution of the equation
\[
\frac{ \mathcal{Z}\left(  \hat{\beta}(  \mu)\right)  }{\mathrm{e}%
^{\hat{\beta}\left(  \mu\right)  \mu}}
=
\frac{\mathbf{E}\left[
\mathrm{e}^{\hat{\beta}\left(  \mu\right)  X}\right]  }{\mathrm{e}%
^{\hat{\beta}\left(  \mu\right)  \mu}}
=
\min_{\beta>0}\frac{\mathbf{E}\left[
\mathrm{e}^{\beta  X}\right]  }{\mathrm{e}%
^{\beta  \mu}}
=
\min_{\beta>0}\frac{\mathcal{Z}(  \beta)  }{\mathrm{e}%
^{\beta  \mu}}.
\]
Let $V\left(
\mu\right)  $ denote the variance of $P^{\mu}.$

Information divergence is given by
\[
D\left(  P^{\mu}\Vert P_{0}\right)  =\int\ln\left(  \frac{\mbox{\textrm{d}%
}P^{\mu}}{\mbox{\textrm{d}}P_{0}}\left(  x\right)  \right)  ~\mbox{\textrm{d}%
}P^{\mu}x.
\]
We see that%
\begin{equation}
\label{*}
D\left(  P^{\mu}\Vert P_{0}\right)
=-\ln\frac{\mathbf{E}\left[
\mathrm{e}^{\hat{\beta}\left(  \mu\right)  X}\right]  }{\mathrm{e}%
^{\hat{\beta}\left(  \mu\right)  \mu}}=\hat{\beta}\left(
\mu\right) \mu-\ln\mathcal{Z}\left(  \hat{\beta}\left(  \mu\right)
\right)  .
\end{equation}

\section{Approximation of tail distributions for non-lattice valued variables}

Introduce the notation
\[
\mu^*:=\sup\{\mu>\mu_{0};
D\left(P^{\mu}\Vert P_{0}\right)<\infty\}=\sup M.
\]
Bahadur and Rao \cite{Bahadur1960a} proved a refined version of
the large deviation bound, but some aspects of their result dates
back to Cram\'er \cite{Cramer1938} and part of it was proved by a
different method by Blackwell and Hodges \cite{Blackwell1959}. For
$\mu^*>\mu>\mu_{0}$,  the Sanov theorem implies that
\[
-\frac{\ln \mathbf{P}\left\{ \frac{1}{n}\sum_{i=1}^{n}X_{i}\geq\mu\right\}}{n}\rightarrow D\left(  P^{\mu}\Vert P_{0}\right)
\mbox{ for
}n\rightarrow\infty.
\]
Bahadur and Rao \cite{Bahadur1960a} verified the
following improvement of the Sanov theorem
\begin{equation}
\mathbf{P}\left\{ \frac{1}{n}\sum_{i=1}^{n}X_{i}\geq\mu\right\}
=
\frac{\exp\left( -nD\left(  P^{\mu}\Vert P_{0}\right)
\right)
}{\left(  {2\pi} nV\left(  \mu\right)  \right)  ^{%
1/2%
}\hat{\beta}\left(  \mu\right)  }
\left( 1+O\left(\frac{1}{\sqrt{n}}\right)\right)\mbox{ for
}n\rightarrow
\infty\label{In}%
\end{equation}
for non lattice random variables.

We will write $D\left( \mu\right)  $ as short for $D\left(  \left.
P^{{\mu}}\right\Vert P_{0}\right)  .$

\begin{theorem}
\label{th2} For $\mu^*>\mu>\mu_{0}$, one has that
\begin{equation}
\mathbf{P}\left\{ \frac{1}{n}\sum_{i=1}^{n}X_{i}\geq\mu\right\}
=
\Phi\left(  -n^{1/2}\left(  2D\left(  \mu-\frac{c_{\mu}}{n}\right)
\right)  ^{1/2}\right)  \left( 1+O\left(\frac{1}{\sqrt{n}}\right)\right)\mbox{
for }n\rightarrow\infty,
\end{equation}
where%
\begin{equation}
c_{\mu}=  \frac{\ln\frac{(2D(  \mu))^{1/2}  }{V\left(  \mu\right)  ^{
1/2}\hat{\beta}\left(  \mu\right)  }}{\hat{\beta}\left(  \mu\right)}.
\label{cce}%
\end{equation}
\end{theorem}

\noindent {\bf Proof.}
The $c_{\mu}$ defined by
(\ref{cce}) satisfies the
equation%
\begin{equation}
\label{fi0}
 \frac{\left(  \frac{2D\left(  \mu\right)  }{V\left(
\mu\right)  }\right)  ^{ 1/2}}{\hat{\beta}\left(  \mu\right)\mathrm{e}^{c_{\mu}\hat{\beta}\left( \mu\right)  }
}=1.
\end{equation}

The tail
probabilities of the
standard Gaussian satisfy%
\[
\frac{\phi\left(  z\right)  }{z}\left(
1-\frac {1}{z^{2}}\right)  \leq\Phi(-z)\leq\frac{\phi\left(
z\right)  }{z}%
\]
for $z>0,$
(cf. Feller {\cite[p. 179]{Feller1957}}),
which implies that
\[
\frac{\exp\left(  -nD\left(
\mu-\frac{c_{\mu}}{n}\right)
\right)  }{\left(  {2\pi} n\right)  ^{1/2}\left(  2D\left(  \mu-\frac{c_{\mu}%
}{n}\right)  \right)  ^{1/2}}
=
\Phi\left(  -n^{1/2}
\left(  2D\left(  \mu-\frac{c_{\mu}}{n}\right)  \right)  ^{1/2}\right)
\left(1+O\left(  \frac{1}{n}\right) \right) ,
\]
and so
\begin{equation}
\label{fi1}
 \frac{\exp\left(  -nD\left(
\mu-\frac{c_{\mu}}{n}\right) \right) }{\left(  {2\pi} n\right)
^{1/2}(  2D(  \mu) ) ^{1/2}}
=
 \Phi\left(  -n^{1/2} \left(
2D\left( \mu-\frac{c_{\mu}}{n}\right) \right)  ^{1/2}\right)
\left(1+O\left(  \frac{1}{n}\right) \right) .
\end{equation}
Because of (\ref{**}) and (\ref{*}), the derivative can be
calculated as
\[
\frac{\mathrm{d}}{\mathrm{d}\mu}D\left(  \mu\right)
=\hat{\beta}\left( \mu\right)  ,
\]
leading to the following Taylor expansion%
\[
D\left(  \mu-\frac{c_{\mu}}{n}\right)  =D\left(  \mu\right)
-\hat{\beta }\left(  \mu\right)  \cdot\frac{c_{\mu}}{n}+O\left(
\frac{1}{n^{2}}\right) .
\]
Thus,
\begin{eqnarray}
 \frac{\exp\left(  -nD\left(
\mu-\frac{c_{\mu}}{n}\right) \right) }{\left(  {2\pi} n\right)
^{1/2}(  2D(  \mu) ) ^{1/2}}
&=&
\frac{\exp\left(  -n\left(D\left(
\mu\right) -\hat{\beta }\left(  \mu\right)
\cdot\frac{c_{\mu}}{n}+O\left( \frac{1}{n^{2}}\right) \right)
\right) }{\left(  {2\pi} n\right) ^{1/2}(  2D(  \mu) ) ^{1/2}}\nonumber\\
&=&
\frac{\exp\left(  -nD(\mu) +\hat{\beta }(  \mu)c_{\mu}+O\left( \frac{1}{n}\right)
\right) }{\left(  {2\pi} n\right) ^{1/2}(  2D(  \mu) )
^{1/2}}\nonumber\\
&=& \frac{\exp\left(  -nD(\mu) \right)
\mathrm{e}^{c_{\mu}\hat{\beta}\left( \mu\right)  }}{\left(  {2\pi}
n\right) ^{1/2}(  2D(  \mu) ) ^{1/2} }\left( 1+O\left(
\frac{1}{n}\right)\right) \label{fi2}
\end{eqnarray}

According to (\ref{In}) we also have
\begin{equation}
\label{fi3} \mathbf{P}\left\{
\frac{1}{n}\sum_{i=1}^{n}X_{i}\geq\mu\right\} = \frac{ \exp\left(
-nD\left( \mu\right) \right)}{\left( {2\pi} nV\left(  \mu\right)
\right)  ^{1/2}
 \hat{\beta}\left(  \mu\right)  }\left(
1+O\left(\frac{1}{\sqrt{n}}\right)\right)\mbox{ for
}n\rightarrow\infty,
\end{equation}
therefore applying (\ref{fi0}), (\ref{fi1}), (\ref{fi2}) and
(\ref{fi3}) the proof of Theorem \ref{th2} is complete. \qed

\begin{remark}
If in the approximation $c_{\mu}$ is replaced by any other constant $c$
then the ratio of the two approximations tends to a number, which
is not equal to $1$:
\begin{eqnarray*}
\frac{\exp\left(  -nD\left( \mu-\frac{c_{\mu}}{n}\right)\right)}
{\exp\left(  -nD\left( \mu-\frac{c}{n}\right)\right)} &=&
\exp\left( -nD\left( \mu-\frac{c_{\mu}}{n}\right)+nD\left(
\mu-\frac{c}{n}\right)\right)\\
&=& \exp\left(\hat{\beta }\left(  \mu\right)  \cdot
(c_{\mu}-c)+O\left(
\frac{1}{n}\right)\right)\\
&\approx&
\exp\left(\hat{\beta }\left(  \mu\right)  \cdot (c_{\mu}-c)\right)\\
&\ne& 1.
\end{eqnarray*}
\end{remark}

\begin{remark}
If $X_{1}$ has a density with respect to the Lebesgue measure then
Bahadur and Rao \cite{Bahadur1960a} proved the stronger result
that
\[
\mathbf{P}\left\{
\frac{1}{n}\sum_{i=1}^{n}X_{i}\geq\mu\right\}
=
\frac{\exp\left(
-nD\left(  P^{\mu}\Vert P_{0}\right)  \right)
}{\left(  {2\pi} nV\left(  \mu\right)  \right)^{%
1/2}\hat{\beta}\left(  \mu\right)  }
\left(1+O\left(  \frac{1}{n}\right)\right).
\]
Using this result we get the following theorem:
If $X_{1}$ has a density with respect to the Lebesgue measure then%
\[
\mathbf{P}\left\{ \frac{1}{n}\sum_{i=1}^{n}X_{i}\geq\mu\right\}
=
\Phi\left(  -n^{1/2}\left(  2D\left(  \mu-\frac{c_{\mu}}{n}\right)
\right)  ^{1/2}\right)  \left( 1+O\left(\frac{1}{n}\right)\right)\mbox{
for }n\rightarrow\infty,
\]
for any $\mu^*>\mu >\mu_{0}$.
\end{remark}

\section{Results for lattice valued variables}

Now assume that $X_{1},X_{2},\dots$ is a sequence of i.i.d. random
variables with values in a lattice of the type $\left\{
kd+\delta\mid k\in \mathbb{Z}\right\}  .$ For such a sequence
Bahadur and Rao \cite{Bahadur1960a} proved that
\begin{equation}
\label{Inn}
\mathbf{P}\left\{
\frac{1}{n}\sum_{i=1}^{n}X_{i}\geq
\mu\right\}
=
\frac{\exp\left(  -nD\left(  \left.  P^{{\mu}%
}\right\Vert P_{0}\right)  \right)  }{\left(  {2\pi} nV\left(
\mu\right)  \right)  ^{1/2}\frac{1-\exp\left(  -d\hat{\beta}\left(  \mu\right)  \right)  }{d}}%
\left(1+O\left(  \frac{1}{n}\right)\right)
\end{equation}
for any  $n $ such that
$\mathbf{P}\left\{
\frac{1}{n}\sum_{i=1}^{n}X_{i}=\mu\right\}  >0.$ We
note that the result (\ref{In}) for non-lattice variables can be
considered as a
limiting version of (\ref{Inn}) for small $d>0$ because%
\[
\frac{1-\exp\left(  -d\beta\right)
}{d}\rightarrow\beta\mbox{ for }d\rightarrow0.
\]

\begin{theorem}
\label{th3} Assume that $X_{1}$ has values in the lattice $\left\{
kd+\delta\mid k\in\mathbb{Z}\right\}  $  and that $\mu^*>\mu> \mu
_{0} $. Then for any
 $n$ such that
$\mathbf{P}\left\{  \frac{1}{n}\sum_{i=1}^{n}X_{i}=\mu\right\}
>0$ one has
\[
\mathbf{P}\left\{ \frac{1}{n}\sum_{i=1}^{n}X_{i}\geq\mu\right\}
=
\Phi\left(  -n^{1/2}\left(  2D\left(  \mu-\frac{c_{\mu}}{n}\right)
\right)  ^{1/2}\right)  \left( 1+O\left(\frac{1}{n}\right)\right)\mbox{
for }n\rightarrow\infty,
\]
where
\[
c_{\mu}=\frac{\ln\frac{(2D(  \mu))^{1/2}  }{V\left(  \mu\right)  ^{1/2}
\frac{1-\exp\left(  -d\hat{\beta}\left(  \mu\right)  \right)  }{d}}%
}{\hat{\beta}\left(  \mu\right)  }.
\]
\end{theorem}

\noindent {\bf Proof.} If $X_{1}$ is lattice valued then the proof
of Theorem \ref{th2} can be modified by replacing
$\hat{\beta}\left(  \mu\right)  $ by $\frac {1-\exp\left(
-d\hat{\beta}\left(  \mu\right)  \right) }{d}$ at the appropriate
places throughout the proof.
There is no modification in the use of a Taylor expansion. \qed

\bigskip

We now turn to the special case, where $X_{1},\dots,X_{n}$
are i.i.d. Bernoulli random variables with
\[
X_{i}=\left\{
\begin{array}
[c]{ll}%
1 & \mbox{with probability }p,\\
0 & \mbox{with probability }1-p.
\end{array}
\right.
\]
In this case $d=1$, and $\sum_{i=1}^{n}X_{i}$ is a binomial $(n,p)$ random variable.
For various refinements of (\ref{Inn}), see Bahadur \cite{Bah60}, Littlewood \cite{Lit69}
and McKay \cite{McK89}.

\begin{cor}
Put
\[
\mu_n:=\lceil n\mu \rceil /n.
\]
Then for $1>\mu>p$ one has that
\[
\mathbf{P}\left\{ \frac{1}{n}\sum_{i=1}^{n}X_{i}\geq\mu\right\}
=
\Phi\left(  -n^{1/2}\left(  2D\left(  \mu_n-\frac{c_{\mu_n}}{n}\right)
\right)  ^{1/2}\right)  \left( 1+O\left(\frac{1}{n}\right)\right)\mbox{
for }n\rightarrow\infty,
\]
where
\[
D(\mu)=D(\mu\Vert p)=\mu\ln\frac{\mu}{p}+(1-\mu)\ln\frac{1-\mu}{1-p}
\]
and
\[
c_{\mu}=
\frac 12 +\frac{\ln\left(  \frac{2D(\mu\Vert p)}{\left(
\mu-p\right) ^{2}}p\left(  1-p\right)  \right)
}{2\ln\frac{\mu\left(  1-p\right) }{p\left(  1-\mu\right)  }}.
\]
\end{cor}

\noindent {\bf Proof.}
Because of the definition of $\mu_n$,
\[
\mathbf{P}\left\{ \frac{1}{n}\sum_{i=1}^{n}X_{i}\geq\mu\right\}
=
\mathbf{P}\left\{ \frac{1}{n}\sum_{i=1}^{n}X_{i}\geq\mu_n\right\},
\]
and the condition $\mathbf{P}\left\{  \frac{1}{n}\sum_{i=1}^{n}X_{i}=\mu_n\right\}
>0$ is satisfied, and so Theorem \ref{th3} implies that
\[
\mathbf{P}\left\{ \frac{1}{n}\sum_{i=1}^{n}X_{i}\geq\mu\right\}
=
\Phi\left(  -n^{1/2}\left(  2D\left(  \mu_n-\frac{c_{\mu_n}}{n}\right)
\right)  ^{1/2}\right)  \left( 1+O\left(\frac{1}{n}\right)\right)\mbox{
for }n\rightarrow\infty.
\]
We have to evaluate $c_{\mu}$.
The distribution $P_{\beta}$ has%
\[
P_{\beta}\left(  X_{i}=1\right)  =\frac{p\mathrm{e}^{\beta}}{1-p+p\mathrm{e}%
^{\beta}}%
\]
which is also the mean of $P_{\beta}$. The equation
\[
\mu=\frac{p\mathrm{e}^{\beta}}{1-p+pe^{\beta}}%
\]
is equivalent to
\[
\mathrm{e}^{\beta}=\frac{\mu\left(  1-p\right)  }{p\left(  1-\mu\right)  }%
\]
implying that%
\[
\frac{1-\mathrm{e}^{-d\beta}}{d}
  =1-\mathrm{e}^{-\beta}
  =1-\frac{p\left(  1-\mu\right)  }{\mu\left(  1-p\right)  }
  =\frac{\mu-p}{\mu\left(  1-p\right)  }.
\]
The variance function is%
\[
V\left(  \mu\right)  =\mu\left(  1-\mu\right)  .
\]
Thus, we have%
\begin{align*}
c_{\mu}  &  =\frac{\ln\left(  \left(  \frac{2D(\mu\Vert
p)}{V\left(
\mu\right)  }\right)  ^{1/2}\frac{1}{1-\mathrm{e}^{  -\hat{\beta}(  \mu) } }\right)  }{\hat{\beta}\left(  \mu\right)  }\\
&  =\frac{\ln\left(  \left(  \frac{2D(\mu\Vert p)}{\mu\left(
1-\mu\right)
}\right)  ^{1/2}\frac{\mu\left(  1-p\right)  }{\mu-p}\right)
}{\ln\frac{\mu\left(
1-p\right)  }{p\left(  1-\mu\right)  }}\\
&  =\frac 12 +\frac{\ln\left(  \frac{2D(\mu\Vert p)}{\left(
\mu-p\right) ^{2}}p\left(  1-p\right)  \right)
}{2\ln\frac{\mu\left(  1-p\right) }{p\left(  1-\mu\right)  }}.
\end{align*}
\qed

\begin{remark}
For $p=1/2$, $0.5<c_{\mu}<0.534$ and Table \ref{tab1} shows some
numerical values for $c_{\mu}\approx 0.5+(\mu-0.5)/12$.

\begin{table}[h]
  \begin{center}\baselineskip=12pt
\begin{tabular}{| c | c | c | c | c | c | c | c |}
\hline
$\mu$ & 0.6 & 0.65 & 0.7 & 0.75 & 0.8 & 0.85 & 0.9\\
\hline
$c_{\mu}$ & 0.508 & 0.512 & 0.516 & 0.520 & 0.524 & 0.528 & 0.532\\
\hline
\end{tabular}
\caption{Numerical values}
\label{tab1}
  \end{center}
\end{table}
\end{remark}

\section{Discussion}

As discussed by Reiczigel, Rejt\H{o} and Tusn\'ady
\cite{Reiczigel2011} and by Harremo\"es and Tusn\'ady
\cite{Harremoes2012} there are some strong indications that these asymptotic
results can be strengthened to sharp inequalities. Such sharp
inequalities would imply the present asymptotic results as
corollaries. We hope that the asymptotics presented here can help
in proving the conjectured sharp inequalities. Related sharp
inequalities have been discussed by Leon and Perron
\cite{Leon2003} and Talagrand \cite{Talagrand1995}. Numerical experiments have also shown that our tail
estimates are useful even for small values of $n.$

\bibliographystyle{elsarticle-num}

\end{document}